\documentclass[12pt]{article}

\usepackage{amsmath,amssymb,amsthm,amscd,a4wide,upref}

\usepackage[enableskew]{youngtab}

\usepackage{hyperref}
\hypersetup{colorlinks,citecolor=blue,filecolor=black,linkcolor=blue,urlcolor=blue}
\usepackage{enumerate}
\usepackage{mathtools}
\usepackage{enumitem}

\usepackage{lmodern}     
\usepackage[T1]{fontenc}

\begin{document}


\newcommand{\ad}{{\rm ad}}
\newcommand{\cri}{{\rm cri}}
\newcommand{\row}{{\rm row}}
\newcommand{\col}{{\rm col}}
\newcommand{\Ann}{{\rm{Ann}\ts}}
\newcommand{\End}{{\rm{End}\ts}}
\newcommand{\Rep}{{\rm{Rep}\ts}}
\newcommand{\Hom}{{\rm{Hom}}}
\newcommand{\Mat}{{\rm{Mat}}}
\newcommand{\ch}{{\rm{ch}\ts}}
\newcommand{\chara}{{\rm{char}\ts}}
\newcommand{\diag}{{\rm diag}}
\newcommand{\st}{{\rm st}}
\newcommand{\non}{\nonumber}
\newcommand{\wt}{\widetilde}
\newcommand{\wh}{\widehat}
\newcommand{\ol}{\overline}
\newcommand{\ot}{\otimes}
\newcommand{\la}{\lambda}
\newcommand{\La}{\Lambda}
\newcommand{\De}{\Delta}
\newcommand{\al}{\alpha}
\newcommand{\be}{\beta}
\newcommand{\ga}{\gamma}
\newcommand{\Ga}{\Gamma}
\newcommand{\ep}{\epsilon}
\newcommand{\ka}{\kappa}
\newcommand{\vk}{\varkappa}
\newcommand{\si}{\sigma}
\newcommand{\vs}{\varsigma}
\newcommand{\vp}{\varphi}
\newcommand{\ta}{\theta}
\newcommand{\de}{\delta}
\newcommand{\ze}{\zeta}
\newcommand{\om}{\omega}
\newcommand{\Om}{\Omega}
\newcommand{\ee}{\epsilon^{}}
\newcommand{\su}{s^{}}
\newcommand{\hra}{\hookrightarrow}
\newcommand{\rar}{\rightarrow}
\newcommand{\lar}{\leftarrow}
\newcommand{\ve}{\varepsilon}
\newcommand{\pr}{^{\tss\prime}}
\newcommand{\ts}{\,}
\newcommand{\vac}{\mathbf{1}}
\newcommand{\vacu}{|0\rangle}
\newcommand{\di}{\partial}
\newcommand{\qin}{q^{-1}}
\newcommand{\tss}{\hspace{1pt}}
\newcommand{\Sr}{ {\rm S}}
\newcommand{\U}{ {\rm U}}
\newcommand{\BL}{ {\overline L}}
\newcommand{\BE}{ {\overline E}}
\newcommand{\BP}{ {\overline P}}
\newcommand{\AAb}{\mathbb{A}\tss}
\newcommand{\CC}{\mathbb{C}\tss}
\newcommand{\KK}{\mathbb{K}\tss}
\newcommand{\QQ}{\mathbb{Q}\tss}
\newcommand{\SSb}{\mathbb{S}\tss}
\newcommand{\TT}{\mathbb{T}\tss}
\newcommand{\ZZ}{\mathbb{Z}\tss}
\newcommand{\DY}{ {\rm DY}}
\newcommand{\X}{ {\rm X}}
\newcommand{\Y}{ {\rm Y}}
\newcommand{\Z}{{\rm Z}}
\newcommand{\ZX}{{\rm ZX}}
\newcommand{\ZY}{{\rm ZY}}
\newcommand{\Ac}{\mathcal{A}}
\newcommand{\Lc}{\mathcal{L}}
\newcommand{\Mc}{\mathcal{M}}
\newcommand{\Pc}{\mathcal{P}}
\newcommand{\Qc}{\mathcal{Q}}
\newcommand{\Rc}{\mathcal{R}}
\newcommand{\Sc}{\mathcal{S}}
\newcommand{\Tc}{\mathcal{T}}
\newcommand{\Bc}{\mathcal{B}}
\newcommand{\Ec}{\mathcal{E}}
\newcommand{\Fc}{\mathcal{F}}
\newcommand{\Gc}{\mathcal{G}}
\newcommand{\Hc}{\mathcal{H}}
\newcommand{\Uc}{\mathcal{U}}
\newcommand{\Vc}{\mathcal{V}}
\newcommand{\Wc}{\mathcal{W}}
\newcommand{\Yc}{\mathcal{Y}}
\newcommand{\Cl}{\mathcal{C}l}
\newcommand{\Ar}{{\rm A}}
\newcommand{\Br}{{\rm B}}
\newcommand{\Ir}{{\rm I}}
\newcommand{\Fr}{{\rm F}}
\newcommand{\Jr}{{\rm J}}
\newcommand{\Or}{{\rm O}}
\newcommand{\GL}{{\rm GL}}
\newcommand{\Spr}{{\rm Sp}}
\newcommand{\Rr}{{\rm R}}
\newcommand{\Zr}{{\rm Z}}
\newcommand{\gl}{\mathfrak{gl}}
\newcommand{\middd}{{\rm mid}}
\newcommand{\ev}{{\rm ev}}
\newcommand{\Pf}{{\rm Pf}}
\newcommand{\Norm}{{\rm Norm\tss}}
\newcommand{\oa}{\mathfrak{o}}
\newcommand{\spa}{\mathfrak{sp}}
\newcommand{\osp}{\mathfrak{osp}}
\newcommand{\f}{\mathfrak{f}}
\newcommand{\g}{\mathfrak{g}}
\newcommand{\h}{\mathfrak h}
\newcommand{\n}{\mathfrak n}
\newcommand{\m}{\mathfrak m}
\newcommand{\se}{\mathfrak{s}}
\newcommand{\z}{\mathfrak{z}}
\newcommand{\Zgot}{\mathfrak{Z}}
\newcommand{\p}{\mathfrak{p}}
\newcommand{\sll}{\mathfrak{sl}}
\newcommand{\agot}{\mathfrak{a}}
\newcommand{\bgot}{\mathfrak{b}}
\newcommand{\qdet}{ {\rm qdet}\ts}
\newcommand{\Ber}{ {\rm Ber}\ts}
\newcommand{\HC}{ {\mathcal HC}}
\newcommand{\cdet}{{\rm cdet}}
\newcommand{\rdet}{{\rm rdet}}
\newcommand{\tr}{ {\rm tr}}
\newcommand{\gr}{ {\rm gr}\ts}
\newcommand{\str}{ {\rm str}}
\newcommand{\loc}{{\rm loc}}
\newcommand{\Gr}{{\rm G}}
\newcommand{\sgn}{ {\rm sgn}\ts}
\newcommand{\sign}{{\rm sgn}}
\newcommand{\ba}{\bar{a}}
\newcommand{\bb}{\bar{b}}
\newcommand{\bi}{\bar{\imath}}
\newcommand{\bj}{\bar{\jmath}}
\newcommand{\bk}{\bar{k}}
\newcommand{\bl}{\bar{l}}
\newcommand{\bp}{\bar{p}}
\newcommand{\hb}{\mathbf{h}}
\newcommand{\Sym}{\mathfrak S}
\newcommand{\fand}{\quad\text{and}\quad}
\newcommand{\Fand}{\qquad\text{and}\qquad}
\newcommand{\For}{\qquad\text{or}\qquad}
\newcommand{\for}{\quad\text{or}\quad}
\newcommand{\grpr}{{\rm gr}^{\tss\prime}\ts}
\newcommand{\degpr}{{\rm deg}^{\tss\prime}\tss}
\newcommand{\bideg}{{\rm bideg}\ts}

\renewcommand{\theequation}{\arabic{section}.\arabic{equation}}

\numberwithin{equation}{section}

\newtheorem{thm}{Theorem}[section]
\newtheorem{lem}[thm]{Lemma}
\newtheorem{prop}[thm]{Proposition}
\newtheorem{cor}[thm]{Corollary}
\newtheorem{conj}[thm]{Conjecture}
\newtheorem*{mthm}{Main Theorem}
\newtheorem*{mthma}{Theorem A}
\newtheorem*{mthmb}{Theorem B}
\newtheorem*{mthmc}{Theorem C}
\newtheorem*{mthmd}{Theorem D}

\theoremstyle{definition}
\newtheorem{defin}[thm]{Definition}

\theoremstyle{remark}
\newtheorem{remark}[thm]{Remark}
\newtheorem{example}[thm]{Example}
\newtheorem{examples}[thm]{Examples}

\newcommand{\bth}{\begin{thm}}
\renewcommand{\eth}{\end{thm}}
\newcommand{\bpr}{\begin{prop}}
\newcommand{\epr}{\end{prop}}
\newcommand{\ble}{\begin{lem}}
\newcommand{\ele}{\end{lem}}
\newcommand{\bco}{\begin{cor}}
\newcommand{\eco}{\end{cor}}
\newcommand{\bde}{\begin{defin}}
\newcommand{\ede}{\end{defin}}
\newcommand{\bex}{\begin{example}}
\newcommand{\eex}{\end{example}}
\newcommand{\bes}{\begin{examples}}
\newcommand{\ees}{\end{examples}}
\newcommand{\bre}{\begin{remark}}
\newcommand{\ere}{\end{remark}}
\newcommand{\bcj}{\begin{conj}}
\newcommand{\ecj}{\end{conj}}

\newcommand{\bal}{\begin{aligned}}
\newcommand{\eal}{\end{aligned}}
\newcommand{\beq}{\begin{equation}}
\newcommand{\eeq}{\end{equation}}
\newcommand{\ben}{\begin{equation*}}
\newcommand{\een}{\end{equation*}}

\newcommand{\bpf}{\begin{proof}}
\newcommand{\epf}{\end{proof}}

\def\beql#1{\begin{equation}\label{#1}}

\newcommand{\Res}{\mathop{\mathrm{Res}}}

\title{\Large\bf Representations of the super-Yangian of type $D(n,m)$}

\author{A. I. Molev}

\date{} 
\maketitle


\begin{abstract}
We consider the classification problem for finite-dimensional irreducible representations of
the Yangians associated with the orthosymplectic Lie superalgebras $\osp_{2n|2m}$ with $n\geqslant 2$.
We give necessary
conditions for an irreducible highest weight representation to be finite-dimensional.
We conjecture that these conditions are also sufficient and prove the conjecture
for a class of representations with linear highest weights. The arguments are based on
a new type of odd reflections for the Yangian associated with $\osp_{2|2}$.
In the Appendix, we construct an isomorphism between the Yangians associated with
the Lie superalgebras $\osp_{2|2}$ and $\gl_{1|2}$.
\end{abstract}



%

\section{Introduction}\label{sec:int}
\setcounter{equation}{0}

In our recent work with E.~Ragoucy~\cite{mr:rb} we gave a conjectural classification
of the finite-dimensional irreducible representations of
the Yangians associated with the orthosymplectic Lie superalgebras $\osp_{2n+1|2m}$.
The conjecture holds in the cases $n=0$ (considered earlier in \cite{m:ry}) and $n=1$
(proved in \cite{mr:rb}). The present paper is concerned with extending the conjecture to the remaining
family of the Yangians associated with the orthosymplectic Lie superalgebras
$\osp_{2n|2m}$ with $n\geqslant 2$
which form the series $D(n,m)$ in the Kac classification of simple Lie superalgebras.
A classification theorem for the Yangians $\Y(\osp_{2|2m})$ was proved earlier in \cite{m:rs}.

By a standard argument, every finite-dimensional irreducible representation of
the Yangian $\Y(\osp_{N|2m})$ with $N=2n$ or $N=2n+1$ is a highest weight representation. Therefore,
the classification problem can be solved by finding necessary
and sufficient conditions on the highest weight for the representation to be
finite-dimensional. The necessity of the conjectural conditions in the case $N=2n+1$
was established in \cite{mr:rb} with the use of the {\em embedding theorem} \cite{m:dt} and
{\em odd reflections} \cite{m:or}. We will use the same instruments in the case $N=2n$
to state a conjectural finite-dimensionality criterion, but they
have to be complemented by a new type of odd reflections. These are certain transformations
of highest weights associated with
the extended Yangian $\Y(\osp_{2|2})$; see Proposition~\ref{prop:odd22}.
By using the combined approaches, we will derive necessary conditions
on the highest weights of $\Y(\osp_{2n|2m})$-modules to be finite-dimensional
(Theorem~\ref{thm:necess}).
We conjecture that these conditions are also sufficient (Conjecture~\ref{conj:fdim}).
We will prove the conjecture
for the {\em linear highest weights} (Theorem~\ref{thm:linear}),
as defined in Sec.~\ref{sec:rl}.

In the Appendix, we consider the Yangian for the general linear Lie superalgebra
$\gl_{1|2}$
and the extended Yangian for the orthosymplectic Lie superalgebra $\osp_{2|2}$
associated with distinguished parity sequences. We produce an explicit isomorphism
between these algebras and use it to give an alternative proof of
the key finite-dimensionality
conditions stated in Corollary~\ref{cor:fdim} for the highest weight representations of
the orthosymplectic Yangian for $\osp_{2|2}$
with the standard parity sequence.

\section{Definitions and preliminaries}
\label{sec:ns}

For given positive integers $m$ and $n$ consider
the {\em parity sequences} $\se=\se_1\dots\se_{m+n}$ of length $m+n$,
where each term $\se_i$ is $0$ or $1$, and the total number of zeros is $n$.
The {\em standard sequence} $\se^{\st}=1\dots 1\tss 0\dots 0$
is defined by $\se_i=1$ for $i=1,\dots,m$
and $\se_i=0$ for $i=m+1,\dots,m+n$.

Suppose a parity sequence $\se$ is fixed. We will simply write $\bi$ to
denote its $i$-th term $\se_i$. Introduce the
involution $i\mapsto i\pr=2n+2m-i+1$ on
the set $\{1,2,\dots,2n+2m\}$ and set $\overline{i'}=\bi$ for $i=1,\dots,m+n$.
Consider the $\ZZ_2$-graded vector space $\CC^{2n|2m}$ over $\CC$ with the basis
$e_1,e_2,\dots,e_{1'}$, where the parity of the basis vector
$e_i$ is defined to be $\bi\mod 2$.
Accordingly, equip
the endomorphism algebra $\End\CC^{2n|2m}$ with the $\ZZ_2$-gradation, where
the parity of the matrix unit $e_{ij}$ is found by
$\bi+\bj\mod 2$.

We will consider {\em even} square matrices with entries in $\ZZ_2$-graded algebras, their
$(i,j)$ entries will have the parity $\bi+\bj\mod 2$.
The algebra of
even matrices over a superalgebra $\Ac$ will be identified with the tensor product algebra
$\End\CC^{2n|2m}\ot\Ac$, so that a matrix $A=[a_{ij}]$ is regarded as the element
\ben
A=\sum_{i,j=1}^{1'}e_{ij}\ot a_{ij}(-1)^{\bi\tss\bj+\bj}\in \End\CC^{2n|2m}\ot\Ac.
\een
We will use the involutive matrix {\em super-transposition} $t$ defined by
$(A^t)_{ij}=a_{j'i'}(-1)^{\bi\bj+\bj}\tss\ta_i\ta_j$,
where for $i=1,\dots,1'$ we set
\ben
\ta_i=\begin{cases} -1\qquad\text{if}\quad i>m+n\fand \bi=1,\\
\phantom{-}1\qquad\text{otherwise.}
\end{cases}
\een

Introduce the permutation operator $P$ by
\begin{align}
P&=\sum_{i,j=1}^{1'} e_{ij}\ot e_{ji}(-1)^{\bj}\in \End\CC^{2n|2m}\ot\End\CC^{2n|2m}
\non\\[-1em]
\intertext{and set}
Q&=\sum_{i,j=1}^{1'} e_{ij}\ot e_{i'j'}(-1)^{\bi\bj}\ts\ta_i\ta_j
\in \End\CC^{2n|2m}\ot\End\CC^{2n|2m}.
\non
\end{align}
The $R$-{\em matrix} associated with $\osp_{2n|2m}$ is the
rational function in $u$ given by
\ben
R(u)=1-\frac{P}{u}+\frac{Q}{u-\ka},\qquad \ka=n-m-1.
\een
This is a super-version of the $R$-matrix
originally found in \cite{zz:rf}.
Following \cite{aacfr:rp}, we
define the {\it extended Yangian\/}
$\X(\osp_{2n|2m})$ corresponding to the standard parity sequence $\se^{\st}$
as a $\ZZ_2$-graded algebra with generators
$t_{ij}^{(r)}$ of parity $\bi+\bj\mod 2$, where $1\leqslant i,j\leqslant 1'$ and $r=1,2,\dots$,
satisfying the following defining relations.
Introduce the formal series
\beql{tiju}
t_{ij}(u)=\de_{ij}+\sum_{r=1}^{\infty}t_{ij}^{(r)}\ts u^{-r}
\in\X(\osp_{2n|2m})[[u^{-1}]]
\eeq
and combine them into the matrix $T(u)=[t_{ij}(u)]$.
Consider the elements of the tensor product algebra
$\End\CC^{2n|2m}\ot\End\CC^{2n|2m}\ot \X(\osp_{2n|2m})[[u^{-1}]]$ given by
\ben
T_1(u)=\sum_{i,j=1}^{1'} e_{ij}\ot 1\ot t_{ij}(u)(-1)^{\bi\tss\bj+\bj}\fand
T_2(u)=\sum_{i,j=1}^{1'} 1\ot e_{ij}\ot t_{ij}(u)(-1)^{\bi\tss\bj+\bj}.
\een
The defining relations for the algebra $\X(\osp_{2n|2m})$ take
the form of the $RTT$-{\em relation}
\beql{RTT}
R(u-v)\ts T_1(u)\ts T_2(v)=T_2(v)\ts T_1(u)\ts R(u-v).
\eeq
The {\it extended Yangian\/} $\X(\osp^{\tss\se}_{2n|2m})$ associated with
an arbitrary parity sequence $\se$ is defined in the same way
(with $\bi$ understood as $\se_i$) and we note that
it is isomorphic to $\X(\osp_{2n|2m})$;
an isomorphism is given by the map
$
t_{ij}(u)\mapsto t_{\si(i),\si(j)}(u)
$
for a suitable permutation $\si$ of the set $\{1,\dots,1'\}$.

The product $T(u-\ka)\ts T^{\tss t}(u)$ is a scalar matrix with
\beql{ttra}
T(u-\ka)\ts T^{\tss t}(u)=c(u)\tss 1,
\eeq
where $c(u)$ is a series in $u^{-1}$. As found in \cite{aacfr:rp},
all the coefficients of the series $c(u)$ belong to
the center $\ZX(\osp_{2n|2m})$ of $\X(\osp_{2n|2m})$. Moreover, they
generate the center, which can be proved by extending the arguments of \cite[Cor.~3.9]{amr:rp}
to the super case as will also appear in \cite{gk:yo}.

The {\em Yangian} $\Y(\osp_{2n|2m})$
is defined as the subalgebra of
$\X(\osp_{2n|2m})$ which
consists of the elements stable under
the automorphisms
\beql{muf}
t_{ij}(u)\mapsto f(u)\ts t_{ij}(u)
\eeq
for all series
$f(u)\in 1+u^{-1}\CC[[u^{-1}]]$.
We have the tensor product decomposition
\beql{tensordecom}
\X(\osp_{2n|2m})=\ZX(\osp_{2n|2m})\ot \Y(\osp_{2n|2m}).
\eeq
The Yangian $\Y(\osp_{2n|2m})$ can also be regarded as the quotient
of $\X(\osp_{2n|2m})$
by the relation $c(u)=1$.

Explicitly, the defining relations \eqref{RTT} can be written
with the use of super-commutator in terms of the series \eqref{tiju} as follows:
\begin{align}
\big[\tss t_{ij}(u),t_{kl}(v)\big]&=\frac{1}{u-v}
\big(t_{kj}(u)\ts t_{il}(v)-t_{kj}(v)\ts t_{il}(u)\big)
(-1)^{\bi\tss\bj+\bi\tss\bk+\bj\tss\bk}
\non\\
{}&-\frac{1}{u-v-\ka}
\Big(\de_{k i\pr}\sum_{p=1}^{1'}\ts t_{pj}(u)\ts t_{p'l}(v)
(-1)^{\bi+\bi\tss\bj+\bj\tss\bp}\ts\ta_i\ta_p
\label{defrel}\\
&\qquad\qquad\quad
{}-\de_{l j\pr}\sum_{p=1}^{1'}\ts t_{k\tss p'}(v)\ts t_{ip}(u)
(-1)^{\bi\tss\bk+\bj\tss\bk+\bi\tss\bp}\ts\ta_{j'}\ta_{p'}\Big).
\non
\end{align}

A version of the Poincar\'e--Birkhoff--Witt theorem holds for the orthosymplectic Yangian
as was pointed out in \cite{aacfr:rp}, while a detailed proof is given in \cite{gk:yo}; cf.
\cite[Sec.~3]{amr:rp}. It states that the algebra $\X(\osp_{2n|2m})$ is generated by
the coefficients of the series $c(u)$ and $t_{ij}(u)$ with the conditions
\ben
\bal
i+j&\leqslant 2n+2m+1\qquad \text{for}\quad i=1,\dots,m,m',\dots,1'\fand\\
i+j&< 2n+2m+1\qquad \text{for}\quad i=m+1,\dots,(m+1)'.
\eal
\een
Moreover, given any total ordering
on the set of the generators, the ordered monomials with the powers of odd generators
not exceeding $1$, form a basis of the algebra.

For a given parity sequence $\se=0\tss\se'$ which begins with $0$ consider
the extended Yangian $\X(\osp^{\tss\se'}_{2n-2|2m})$ and
for the parity sequence $\se=1\tss\se'$ beginning with $1$ consider
the extended Yangian $\X(\osp^{\tss\se'}_{2n|2m-2})$. In both cases, let
the indices
of the generators $t_{ij}^{(r)}$ of these algebras range over the sets
$2\leqslant i,j\leqslant 2\pr$ and $r=1,2,\dots$. The following
embedding properties were proved in \cite[Thm~3.1]{m:dt} for a standard
parity sequence, and the proof extends to arbitrary sequences $\se$ with only minor
modifications which are listed in \cite[Sec.~3.4]{ft:oy}.
The mapping
\beql{embedgen}
t_{ij}(u)\mapsto t_{ij}(u)-t_{i1}(u)\ts t_{11}(u)^{-1}\tss t_{1j}(u),\qquad 2\leqslant i,j\leqslant 2\pr,
\eeq
defines injective
homomorphisms
\beql{emb}
\X(\osp^{\tss\se'}_{2n-2|2m})\hra \X(\osp^{\tss\se}_{2n|2m})
\Fand
\X(\osp^{\tss\se'}_{2n|2m-2})\hra \X(\osp^{\tss\se}_{2n|2m})
\eeq
in the cases
$\se=0\tss\se'$ and $\se=1\tss\se'$, respectively.

The extended Yangian $\X(\osp_{2n|2m})$ is a Hopf algebra with the coproduct
defined by
\beql{Delta}
\De: t_{ij}(u)\mapsto \sum_{k=1}^{1'} t_{ik}(u)\ot t_{kj}(u).
\eeq
The image of the series $c(u)$ is found by the relation
$\De:c(u)\mapsto c(u)\ot c(u)$ and so the Yangian
$\Y(\osp_{2n|2m})$ inherits the Hopf algebra structure from $\X(\osp_{2n|2m})$.

\section{Highest weight representations}
\label{sec:hw}

We will keep a parity sequence $\se$ fixed.
A representation $V$ of the algebra $\X(\osp^{\tss\se}_{2n|2m})$
is called a {\em highest weight representation}
if there exists a nonzero vector
$\xi\in V$ such that $V$ is generated by $\xi$,
\begin{alignat}{2}
t_{ij}(u)\ts\xi&=0 \qquad &&\text{for}
\quad 1\leqslant i<j\leqslant 1', \qquad \text{and}\non\\
t_{ii}(u)\ts\xi&=\la_i(u)\ts\xi \qquad &&\text{for}
\quad i=1,\dots,1',
\label{trianb}
\end{alignat}
for some formal series
\beql{laiu}
\la_i(u)\in 1+u^{-1}\CC[[u^{-1}]].
\eeq
The vector $\xi$ is called the {\em highest vector}
of $V$.

\bpr\label{prop:nontrvm}
The series $\la_i(u)$ associated with a highest weight representation $V$
satisfy
the consistency conditions
\begin{multline}\label{nontrvm}
\la_i(u)\tss \la_{i\pr}\big(u-n+m+1+(-1)^{\bar 1}+\dots+(-1)^{\bi}\big)\\
{}=\la_{i+1}(u)\tss \la_{(i+1)'}\big(u-n+m+1+(-1)^{\bar 1}+\dots+(-1)^{\bi}\big)
\end{multline}
for $i=1,\dots,m+n-1$.
Moreover, the coefficients of the series $c(u)$ act in the representation
$V$ as the multiplications by scalars
determined by
$
c(u)\mapsto \la_1(u)\tss \la_{1'}(u-n+m+1).
$
\epr

\bpf
Relation \eqref{nontrvm} for $i=1$ follows by applying \eqref{defrel} to calculate
$t_{12}(u)\ts t_{1'2'}(v)\ts \xi$ and then setting
$v=u-\ka+(-1)^{\bar 1}$; cf. \cite[Prop.~4.4]{m:ry}. The formula
for the remaining values of $i$ is then implied by the embedding properties \eqref{emb}.
The eigenvalues of the coefficients of the series
$c(u)$ are found by taking the $(1',1')$ entry in
the matrix relation \eqref{ttra}.
\epf

As Proposition~\ref{prop:nontrvm} shows, the series $\la_i(u)$ in \eqref{trianb}
with $i>m+n+1$ are uniquely
determined by the first $m+n+1$ series. We will call the corresponding tuple
$\la(u)=(\la_{1}(u),\dots,\la_{m+n+1}(u))$
the {\em highest weight\/} of $V$.

Given an arbitrary tuple $\la(u)=(\la_{1}(u),\dots,\la_{m+n+1}(u))$
of formal series of the form \eqref{laiu}, define
the {\em Verma module} $M(\la(u))$ as the quotient of the algebra $\X(\osp^{\tss\se}_{2n|2m})$ by
the left ideal generated by all coefficients of the series $t_{ij}(u)$
with $1\leqslant i<j\leqslant 1'$ and $t_{ii}(u)-\la_i(u)$ for
$i=1,\dots,m+n+1$. The Poincar\'e--Birkhoff--Witt theorem for the algebra $\X(\osp^{\tss\se}_{2n|2m})$
implies that the Verma module $M(\la(u))$
is nonzero, and we denote by $L(\la(u))$ its irreducible quotient.
It is clear that the isomorphism
class of $L(\la(u))$ is determined by $\la(u)$.

\bpr\label{prop:fdhw}
Every finite-dimensional irreducible representation of the extended
Yangian $\X(\osp^{\tss\se}_{2n|2m})$
is isomorphic to the highest weight representation $L(\la(u))$ for a certain highest weight
$\la(u)=(\la_{1}(u),\dots,\la_{m+n+1}(u))$.
\epr

\bpf
The proof follows by a standard argument;
cf. \cite[Thm~5.1]{amr:rp}.
\epf

Recall from \cite{n:qb} that the {\em Yangian} $\Y(\gl_{n|m})$
associated with the general linear Lie superalgebra
$\gl_{n|m}$
(corresponding to the standard parity sequence $\se^{\st}$) is defined as the $\ZZ_2$-graded algebra
with generators $\bar t_{ij}^{\ts(r)}$
of parity $\bi+\bj\mod 2$, where $1\leqslant i,j\leqslant m+n$ and $r=1,2,\dots$,
while
\ben
\bi=\begin{cases} 1\qquad\text{for}\quad i=1,\dots,m,\\
0\qquad\text{for}\quad i=m+1,\dots,m+n.
\end{cases}
\een
The defining relations can be written in terms of the generating series
\beql{btexpa}
\bar t_{ij}(u)=\de_{ij}+\sum_{r=1}^{\infty}\bar t_{ij}^{\ts(r)}\ts u^{-r}
\in\Y(\gl_{n|m})[[u^{-1}]]
\eeq
and they have the form
\beql{gldefrel}
\big[\tss \bar t_{ij}(u),\bar t_{kl}(v)\big]=\frac{1}{u-v}
\big(\bar t_{kj}(u)\ts \bar t_{il}(v)-\bar t_{kj}(v)\ts \bar t_{il}(u)\big)
(-1)^{\bi\tss\bj+\bi\tss\bk+\bj\tss\bk}.
\eeq
By the Poincar\'e--Birkhoff--Witt theorem for the Yangian $\Y(\gl_{n|m})$ \cite{g:gd} and
the orthosymplectic Yangian (Sec.~\ref{sec:ns}), we may regard $\Y(\gl_{n|m})$
as a subalgebra of $\X(\osp_{2n|2m})$
via the embedding
\beql{embedya}
\Y(\gl_{n|m})\hra \X(\osp_{2n|2m}),
\eeq
given by
\ben
\bar t_{ij}(u)\mapsto t_{ij}(u),\qquad 1\leqslant i,j\leqslant m+n.
\een

The following simple observation will be used below:
if a highest weight representation $L(\la(u))$ of $\X(\osp_{2n|2m})$ is finite-dimensional, then so is
the $\Y(\gl_{n|m})$-module $\Y(\gl_{n|m})\xi$ defined via the embedding \eqref{embedya}.
Hence the tuple
\beql{glnmt}
\big(\la_1(u),\dots,\la_{m}(u),\la_{m+1}(u),\dots,\la_{m+n}(u)\big)
\eeq
satisfies the finite-dimensionality conditions obtained in \cite{zh:sy}.

\subsection{Odd reflections in the extended Yangian $\X(\osp_{2|2})$}
\label{subsec:or}

We will need a correspondence between the highest weights of the representations
of the extended Yangian $\X(\osp_{2|2})$ associated
with different parity sequences. It will complement the correspondence
between the highest weights of the representations
of the Yangian $\Y(\gl_{n|m})$ obtained earlier in \cite{m:or}.
This correspondence was used therein to give finite-dimensionality conditions
for highest weight representations of $\Y(\gl_{n|m})$ associated
with arbitrary parity sequences. This generalized the results of \cite{zh:sy},
where the standard parity sequence was considered.

We have two possible parity sequences $01$ and $10$ for $\osp_{2|2}$,
and the isomorphism
\beql{isom22}
\vs:\X(\osp^{01}_{2|2})\to\X(\osp^{10}_{2|2})
\eeq
given by the mapping
\ben
t_{ij}(u)\mapsto t_{\si(i),\si(j)}(u),\qquad i,j\in\{1,2,2^{\tss\prime},1'\},
\een
for $\si=(12)(2^{\tss\prime}1')$ (note that here $2\pr=3$ and $1'=4$).
The mapping is illustrated by the diagram, where the $(i,j)$ box represents
the generator $t_{ij}(u)$:

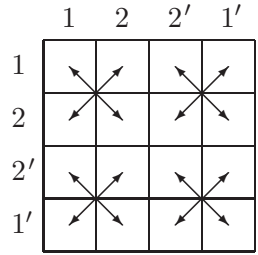
\bigskip

\begin{center}
\begin{picture}(80,80)
\thinlines

\put(0,0){\line(0,1){80}}
\put(0,0){\line(1,0){80}}
\put(0,20){\line(1,0){80}}
\put(20,0){\line(0,1){80}}
\put(0,40){\line(1,0){80}}
\put(40,0){\line(0,1){80}}
\put(0,60){\line(1,0){80}}
\put(60,0){\line(0,1){80}}
\put(0,80){\line(1,0){80}}
\put(80,0){\line(0,1){80}}

\put(10,10){\vector(1,1){20}}
\put(30,30){\vector(-1,-1){20}}
\put(10,30){\vector(1,-1){20}}
\put(30,10){\vector(-1,1){20}}

\put(10,50){\vector(1,1){20}}
\put(30,70){\vector(-1,-1){20}}
\put(10,70){\vector(1,-1){20}}
\put(30,50){\vector(-1,1){20}}

\put(50,10){\vector(1,1){20}}
\put(70,30){\vector(-1,-1){20}}
\put(50,30){\vector(1,-1){20}}
\put(70,10){\vector(-1,1){20}}

\put(50,50){\vector(1,1){20}}
\put(70,70){\vector(-1,-1){20}}
\put(50,70){\vector(1,-1){20}}
\put(70,50){\vector(-1,1){20}}

\put(7,85){\small$1$}
\put(27,85){\small$2$}
\put(47,85){\small$2^{\tss\prime}$}
\put(67,85){\small$1'$}

\put(-12,7){\small$1'$}
\put(-12,27){\small$2^{\tss\prime}$}
\put(-12,47){\small$2$}
\put(-12,67){\small$1$}

\end{picture}
\end{center}

Given any $\X(\osp^{10}_{2|2})$-module $V$, denote by $V^{\vs}$ the $\X(\osp^{01}_{2|2})$-module
obtained by the composition of the action of the algebra $\X(\osp^{10}_{2|2})$ on $V$
with the isomorphism \eqref{isom22}.

The following is a version of odd reflections for the extended Yangian $\X(\osp_{2|2})$.

\bpr\label{prop:odd22}
Suppose that the first two components $\la_1(u)$ and $\la_2(u)$ of the highest weight
$\la(u)=(\la_{1}(u),\la_{2}(u),\la_{2'}(u))$ of the $\X(\osp^{10}_{2|2})$-module $L(\la(u))$
have the form
\ben
\la_1(u)=(1+\al_1u^{-1})\dots(1+\al_pu^{-1})\fand
\la_2(u)=(1+\be_1u^{-1})\dots(1+\be_pu^{-1})
\een
for some nonnegative integer $p$ with the conditions $\al_i\ne\be_j$ for all
$i,j\in\{1,\dots,p\}$. Then the $\X(\osp^{01}_{2|2})$-module $L(\la(u))^{\vs}$
is isomorphic to the highest weight module $L(\wt\la(u))$, where the components
of $\wt\la(u)$ are given by
\beql{laot}
\wt\la_1(u)=\Big(\frac{u+1}{u}\Big)^p\ts\la_2(u+1),\qquad
\wt\la_2(u)=\Big(\frac{u+1}{u}\Big)^p\ts\la_1(u+1)
\eeq
and
\beql{latprim}
\wt\la_{2'}(u)=\Big(\frac{u-1}{u}\Big)^p\ts\ts\frac{\la_{2}(u-1)\tss\la_{2'}(u)}{\la_1(u)}.
\eeq
\epr

\bpf
Identify the Yangian $\Y(\gl_{1|1})$ associated with the standard parity sequence
$\se^{\st}=10$ with the subalgebra of $\X(\osp_{2|2})$
via the embedding \eqref{embedya} for $m=n=1$. The cyclic span $\BL:=\Y(\gl_{1|1})\xi$
of the highest vector $\xi\in L(\la(u))$ is a highest weight module
over $\Y(\gl_{1|1})$ with the highest weight $(\la_1(u),\la_2(u))$. Moreover,
we can also observe that all generators $t_{21}^{(r)}$ with $r>p$ act as zero
operators in $\BL$. This follows from
the relation $t_{21}^{(r)}\xi=0$ in $L(\la(u))$ for these values of $r$,
which holds since
the vector $t_{21}^{(r)}\xi\in M(\la(u))$ is annihilated by
all operators $t_{ij}(u)$ with $1\leqslant i<j\leqslant 1'$.

Since the operator $T_{21}(u):=u^p\tss t_{21}(u)$ in $\BL$ is a polynomial in $u$,
we can introduce the vector $\ze\in\BL$ by
\ben
\ze=T_{21}(-\al_1)\dots T_{21}(-\al_p)\ts\xi.
\een
According to \cite[Prop.~4.1]{m:or} (applied to the opposite parity sequence), this vector
is nonzero and the following relations hold:
\begin{align}
t_{11}(u)\ts\ze&=\Big(\frac{u+1}{u}\Big)^p\ts\la_1(u+1)\ts\ze,
\label{too}\\[0.4em]
t_{22}(u)\ts\ze&=\Big(\frac{u+1}{u}\Big)^p\ts\la_2(u+1)\ts\ze,
\label{ttt}\\[0.4em]
t_{21}(u)\ts\ze&=0.
\label{tto}
\end{align}
To find the action of other generators of $\X(\osp_{2|2})$ on the vector $\ze$, use
the quasideterminant formulas for the Drinfeld generators
introduced in \cite[Sec.~4]{m:dt}:
\beql{gauh}
h_1(u)=t_{11}(u),\quad
h_2(u)=\begin{vmatrix} t_{1\tss 1}(u)&t_{1\tss 2}(u)\\
                         t_{2\tss 1}(u)&\boxed{t_{2\tss 2}(u)}
           \end{vmatrix},
\quad
h_{3}(u)=\begin{vmatrix} t_{1\tss 1}(u)&t_{1\ts 2}(u)&t_{1\tss 3}(u)\\
                         t_{2\tss 1}(u)&t_{2\tss 2}(u)&t_{2\tss 3}(u)\\
                         t_{31}(u)&t_{32}(u)&\boxed{t_{33}(u)}
           \end{vmatrix}.
\eeq
According to \cite[Thm~5.3]{m:dt}, the series $c(u)$ defined
in \eqref{ttra} is calculated by
\ben
c(u)=\frac{h_1(u)}{h_1(u+1)}\ts h_2(u+1)\tss h_{3}(u+1).
\een
On the other hand, the {\em quantum Berezinian} of \cite{n:qb}
for the subalgebra isomorphic to $\Y(\gl_{1|1})$
is found by the formula $b(u)=h_1(u)\tss h_2(u)^{-1}$;
see \cite[Sec.~7]{g:gd}.
All coefficients of $b(u)$ belong to the center of the subalgebra.
Therefore, the series $h_{3}(u)$ can be written as
\beql{htwopr}
h_{3}(u)=c(u-1)\tss b(u)\tss h_1(u-1)^{-1}.
\eeq
Since $b(u)$ commutes with $t_{21}(v)$, and
\ben
h_2(u)=t_{22}(u)-t_{21}(u)\ts t_{11}(u)^{-1}\tss t_{12}(u),
\een
for the action on the vector $\ze$ we have
\ben
b(u)\tss\ze=\la_1(u)\tss \la_2(u)^{-1}\tss\ze,
\een
while Proposition~\ref{prop:nontrvm} gives
\ben
c(u)\tss \ze=\la_1(u)\tss \la_{1'}(u+1)\tss \ze.
\een
Hence, \eqref{too} and \eqref{htwopr} yield
\beql{htpri}
h_{3}(u)\tss\ze=\Big(\frac{u-1}{u}\Big)^p\ts\ts\frac{\la_{1}(u-1)\tss\la_{1'}(u)}{\la_2(u)}\ts\ze
=\Big(\frac{u-1}{u}\Big)^p\ts\ts\frac{\la_{1}(u-1)\tss\la_{2'}(u)}{\la_1(u)}\ts\ze,
\eeq
where the second equality follows from the particular case $i=1$ of \eqref{nontrvm}
which reads
\beql{zeoo}
\la_1(u)\la_{1'}(u)=\la_2(u)\la_{2'}(u).
\eeq

The defining relations \eqref{defrel} imply that $t_{21}(u)\ts t_{21}(u+1)=0$.
Therefore, the vector $\ze$ belongs to the subspace $L^{\circ}$ of $L(\la(u))$
spanned by the vectors
\beql{vecq}
t_{21}^{(r_1)}\dots t_{21}^{(r_q)}\ts\xi,\qquad 1\leqslant r_1<\dots<r_q\leqslant p,
\eeq
where $q=1,\dots,p$.

\ble\label{lem:loanni}
The subspace $L^{\circ}$ is annihilated by the action of all coefficients of the series
$t_{12'}(u)$, $t_{21'}(u)$, $t_{11'}(u)$ and $t_{22'}(u)$.
\ele

\bpf
The defining relations \eqref{defrel} give
for $i=1,2$:
\ben
\big[t_{i\tss 2'}(u),t_{21}(v)\big]=(-1)^{\bi+1}\ts\big[t_{21}(v),t_{i\tss 2'}(u)\big]
=-\tss\frac{1}{u-v}
\big(t_{i\tss 1}(u)\ts t_{2\tss 2'}(v)-t_{i\tss 1}(v)\ts t_{2\tss 2'}(u)\big).
\een
Now an easy induction on $q$ shows that the vector \eqref{vecq}
is annihilated by $t_{1\tss 2'}(u)$ and $t_{2\tss 2'}(u)$. Similarly,
by the defining relations we have
\ben
\big[t_{1\tss 1}(u),t_{21}(v)\big]=-\frac{1}{u-v}
\big(t_{2\tss 1}(u)\ts t_{1\tss 1}(v)-t_{2\tss 1}(v)\ts t_{1\tss 1}(u)\big)
\een
and
\ben
\big[t_{2\tss 1}(v),t_{22}(u)\big]=\frac{1}{u-v}
\big(t_{2\tss 1}(u)\ts t_{2\tss 2}(v)-t_{2\tss 1}(v)\ts t_{2\tss 2}(u)\big).
\een
It follows by induction on $q$ in \eqref{vecq}
that the subspace $L^{\circ}\subset L(\la(u))$ is invariant
under the action of $t_{1\tss 1}(u)$ and $t_{2\tss 2}(u)$.
As a consequence of \eqref{ttra}, we get the symmetry
relation $t^{(1)}_{2\tss 1'}=-t^{(1)}_{1\tss 2'}$. So
we derive from the defining relations
that
\ben
t_{21'}(u)=-\big[t^{(1)}_{2\tss 1'},t_{22}(u)\big]=\big[t^{(1)}_{1\tss 2'},t_{22}(u)\big].
\een
Hence $t_{21'}(u)\ts L^{\circ}=0$.
Finally, the property $t_{11'}(u)\ts L^{\circ}=0$ follows from
\ben
t_{11'}(u)=2\ts\big[t^{(1)}_{1\tss 1'},t_{11}(u)\big]\Fand
\big[t^{(1)}_{1\tss 1'},t_{21}(u)\big]=2\ts t_{21'}(u)
\een
and yet another use of induction on $q$ in \eqref{vecq}.
\epf

Returning to the proof of the proposition, observe that
the quasideterminant formula for $h_{3}(u)$ can be written in the form
\ben
h_{3}(u)=t_{33}(u)-\sum_{i,j=1}^2 t_{3i}(u)\tss A_{ij}(u)\ts t_{j\tss 3}(u),
\een
where $A_{ij}(u)$ denotes the $(i,j)$-entry of the matrix
\beql{htta}
A(u)=\begin{bmatrix} t_{1\tss 1}(u)&t_{1\tss 2}(u)\\
                         t_{2\tss 1}(u)&t_{2\tss 2}(u)
           \end{bmatrix}^{-1}.
\eeq
Hence \eqref{htpri} and Lemma~\ref{lem:loanni} imply
\ben
t_{2'2'}(u)\ts\ze =\Big(\frac{u-1}{u}\Big)^p\ts\ts
\frac{\la_{1}(u-1)\tss\la_{2'}(u)}{\la_1(u)}\ts\ze.
\een

Thus, taking into account relations \eqref{too}--\eqref{tto}, we can now conclude that
the $\X(\osp^{01}_{2|2})$-module $L(\la(u))^{\vs}$ is a highest weight module
with the highest vector $\ze$. Moreover, the components of the highest weight
$\wt\la(u)$ are given by
\eqref{laot} together with
\ben
\wt\la_{1'}(u)=\Big(\frac{u-1}{u}\Big)^p\ts\ts\frac{\la_{1}(u-1)\tss\la_{1'}(u)}{\la_2(u)}.
\een
The remaining component $\wt\la_{2'}(u)$ is given by \eqref{latprim}, as follows from the
relation
\ben
\wt\la_{1}(u)\wt\la_{1'}(u+2)=\wt\la_{2}(u)\wt\la_{2'}(u+2),
\een
which is a particular case of \eqref{nontrvm}.
\epf

Recalling that $\X(\osp_{2|2})=\X(\osp^{10}_{2|2})$,
consider the $\X(\osp_{2|2})$-module $L(\la(u))$ whose highest weight $\la(u)$ satisfies
the assumptions of Proposition~\ref{prop:odd22}.

\bco\label{cor:fdim}
The representation $L(\la(u))$ is finite-dimensional if and only if there exists a monic
polynomial $P(u)$ in $u$ such that
\beql{fdco}
\Big(\frac{u-1}{u+1}\Big)^p\ts\ts\frac{\la_{2}(u-1)\tss\la_{2'}(u)}{\la_1(u)\tss \la_1(u+1)}
=\frac{P(u+2)}{P(u)}.
\eeq
\eco

\bpf
By Proposition~\ref{prop:odd22}, the $\X(\osp_{2|2})$-module $L(\la(u))$
is finite-dimensional if and only if the $\X(\osp^{01}_{2|2})$-module $L(\wt\la(u))$
is finite-dimensional. According to the Main Theorem of \cite{m:rs},
this holds if and only if
\ben
\frac{\wt\la_{2'}(u)}{\wt\la_2(u)}
=\frac{P(u+2)}{P(u)}
\een
for a monic polynomial $P(u)$. By \eqref{laot}
and \eqref{latprim}, this is equivalent to \eqref{fdco}.
\epf

\bre\label{rem:fdcom}\quad (i)
Note that if an arbitrary highest weight $\X(\osp_{2|2})$-module $L(\la(u))$
is finite-dimensional, then so is the cyclic span $\Y(\gl_{1|1})\ts\xi$
of the highest vector $\xi$ with respect to the subalgebra $\Y(\gl_{1|1})$.
The cyclic span is a highest weight module over $\Y(\gl_{1|1})$
with the highest weight $(\la_1(u),\la_2(u))$. Hence the ratio $\la_1(u)/\la_2(u)$
is a rational function in $u^{-1}$ due to \cite{zh:rs}. Therefore,
by taking the composition of the $\X(\osp_{2|2})$-module $L(\la(u))$
with a suitable automorphism of the extended Yangian of the form \eqref{muf},
we can get a finite-dimensional module whose highest weight satisfies
the assumptions of Proposition~\ref{prop:odd22}.
Hence, Corollary~\ref{cor:fdim} provides a finite-dimensionality criterion
for an arbitrary highest weight module $L(\la(u))$ over $\X(\osp_{2|2})$.

(ii) Due to the embedding $\X(\oa_2)\hra\X(\osp_{2|2})$ defined in
\eqref{embedgen} and the relations $t_{22'}(u)=t_{2'2}(u)=0$ in $\X(\oa_2)$, we have
the identities
\ben
t_{22'}(u)=t_{21}(u)\ts t_{11}(u)^{-1}\tss t_{12'}(u)\Fand
t_{2'2}(u)=t_{2'1}(u)\ts t_{11}(u)^{-1}\tss t_{12}(u)
\een
in $\X(\osp_{2|2})$. On the other hand,
the mapping $t_{ij}(u)\mapsto t_{\si(i),\si(j)}(u)$
with $\si=(2\ts 2')$ defines an automorphism of the extended Yangian $\X(\osp_{2|2})$.
Hence, it follows from these identities that the composition
of the $\X(\osp_{2|2})$-module $L(\la_1(u),\la_2(u),\la_{2'}(u))$ with the
automorphism is isomorphic to the module $L(\la_1(u),\la_{2'}(u),\la_2(u))$.
This implies that the condition \eqref{fdco} is symmetric with respect to
the transposition $\la_2(u)\leftrightarrow\la_{2'}(u)$ (although the value of the
parameter $p$ can be changed).
We will use the notation
\beql{fdosp}
\la_2(u)\rightrightarrows\la_1(u)\leftleftarrows\la_{2'}(u)
\eeq
to indicate that the finite-dimensionality conditions for the $\X(\osp_{2|2})$-module $L(\la(u))$
are satisfied by the highest weight $\la(u)=(\la_1(u),\la_2(u),\la_{2'}(u))$.
\qed
\ere

\subsection{Necessary conditions}
\label{subsec:thmnec}

Consider the irreducible highest representation $L(\la(u))$ of
the extended Yangian $\X(\osp_{2n|2m})$ with $n\geqslant 2$,
whose highest weight is an arbitrary tuple of formal series
\beql{lau}
\la(u)=\big(\la_1(u),\dots,\la_{m}(u),\la_{m+1}(u),\dots,\la_{m+n}(u),\la_{(m+n)'}(u)\big).
\eeq

As we pointed out above, if $\dim L(\la(u))<\infty$, then
the tuple \eqref{glnmt} satisfies the conditions of \cite{zh:sy}.
In particular, the ratio $\la_m(u)/\la_{m+1}(u)$ is a rational function in $u^{-1}$.
This makes it possible to apply
the $A$ type
odd reflections as introduced in \cite{m:or}.
These are transformations which apply to
pairs $(\al(u),\be(u))$ of formal series of the form
\ben
\bal
\al(u)&=(1+\al_1u^{-1})\dots (1+\al_pu^{-1})\ts\ga(u),\\
\be(u)&=(1+\be_1u^{-1})\dots (1+\be_p\tss u^{-1})\ts\ga(u),
\eal
\een
where $\al_i\ne\be_j$ for all $i,j$, and $\ga(u)\in 1+u^{-1}\CC[[u^{-1}]]$.
The odd reflection is the transformation
\beql{oddrefl}
\big(\al(u),\be(u)\big)\mapsto \big(\be^{[1]}(u),\al^{[1]}(u)\big),
\eeq
where
\ben
\bal
\al^{[1]}(u)&=\big(1+(\al_1+1)\tss u^{-1}\big)\dots \big(1+(\al_p+1)\tss u^{-1}\big)\ts\ga(u),\\
\be^{[1]}(u)&=\big(1+(\be_1+1)\tss u^{-1}\big)\dots \big(1+(\be_p+1)\tss u^{-1}\big)\ts\ga(u).
\eal
\een

An additional series $\la^{[n-1]}_{m}(u)$ for the tuple \eqref{lau}
is obtained by applying the sequence of odd reflections \eqref{oddrefl}
beginning with
\beql{lammpo}
\big(\la_m(u),\la_{m+1}(u)\big)\mapsto \big(\la^{[1]}_{m+1}(u),\la^{[1]}_{m}(u)\big),
\eeq
then using the values of $i=1,\dots,n-2$ and applying \eqref{oddrefl}
consecutively
to determine new series $\la^{[2]}_{m}(u),\dots,\la^{[n-1]}_{m}(u)$ by
\beql{oddlami}
\big(\la^{[i]}_{m}(u),\la_{m+i+1}(u)\big)\mapsto \big(\la^{[1]}_{m+i+1}(u),\la^{[i+1]}_{m}(u)\big).
\eeq

In addition to \eqref{fdosp}, we will also use the notation $\la(u)\rar \mu(u)$ or
$\mu(u)\lar \la(u)$ for two formal power series in $u^{-1}$ to indicate that
\ben
\frac{\la(u)}{\mu(u)}=\frac{Q(u+1)}{Q(u)}
\een
for a monic polynomial $Q(u)$ in $u$.

We can now state necessary conditions for finite-dimensionality of the highest weight representations
of the extended Yangian $\X(\osp_{2n|2m})$ with $n\geqslant 2$.

\bth\label{thm:necess}
If the representation $L(\la(u))$ with the highest weight of the form \eqref{lau}
is finite-dimensional then $\la_m(u)/\la_{m+1}(u)$ is a rational function in $u^{-1}$
and the following conditions hold:
\beql{condgl}
\la_1(u)\lar\la_2(u)\lar\cdots\lar\la_m(u)
\eeq
together with
\beql{condo}
\substack{{\displaystyle{\la_{m+1}(u)\rar\cdots\rar\la_{m+n-1}(u)\rar\la_{m+n}(u)}}
\\{\qquad\displaystyle\downarrow}
\\\substack{\qquad\qquad\displaystyle{\la_{(m+n)'{\displaystyle(u)}}}}}
\eeq
and
\beql{condosp}
\la_{m+n}(u)\ \rightrightarrows\ \la^{[n-1]}_{m}(u)\ \leftleftarrows\ \la_{(m+n)'}(u).
\eeq
\eth

\bpf
First we use the embedding \eqref{embedya} to conclude by \cite{zh:sy} that
$\la_m(u)/\la_{m+1}(u)$ is a rational function in $u^{-1}$, while
the tuple \eqref{glnmt} satisfies both conditions \eqref{condgl}
and the top line of conditions in \eqref{condo}. Similarly, the remaining condition
in \eqref{condo} follows by using the embedding
\beql{embedyad}
\X(\oa_{2n})\hra\X(\osp_{2n|2m})
\eeq
given in \cite[Cor.~3.2]{m:dt}. Indeed, the $\X(\oa_{2n})$-module $\X(\oa_{2n})\ts\xi$
is a finite-dimensional highest weight
representation with the highest weight
$(\la_{m+1}(u),\dots,\la_{m+n}(u),\la_{(m+n)'}(u))$ so that the condition
is implied by the finite-dimensionality criterion of \cite[Thm~5.16]{amr:rp}.

Conditions \eqref{condgl} and \eqref{condo} imply that
by twisting
the representation $L(\la(u))$ by a suitable automorphism \eqref{muf},
if necessary, we may assume that all components $\la_i(u)$ of
the highest weight $\la(u)$ are polynomials in $u^{-1}$; that is,
for some $p\in\ZZ_+$ we have
\beql{hwdec}
\la_i(u)=(1+\la^{(1)}_i u^{-1})\dots (1+\la^{(p)}_i u^{-1}),\qquad i=1,\dots,m+n+1,
\eeq
with $\la^{(r)}_i\in\CC$ for $r=1,\dots,p$.

Consider the parity sequence $\se$ obtained from $\se^{\st}$ by replacing
the subsequence $\se_m\se_{m+1}=10$ with $01$. To calculate
the highest weight of the module $L(\la(u))$ associated with $\se$,
apply the corresponding odd reflection by using \cite[Thm~4.4]{m:or}.
We will assume that the parameters $\la^{(r)}_m$ and $\la^{(r)}_{m+1}$
are numbered in such a way that
\beql{laeq}
\la_m^{(r)}=\la_{m+1}^{(r)}\qquad \text{for all}\quad r=k+1,\dots,p
\eeq
for certain $k\in\{0,1,\dots,p\}$ (relation \eqref{laeq} is vacuous for $k=p$),
while $\la_m^{(r)}\ne\la_{m+1}^{(s)}$
for all $1\leqslant r,s\leqslant k$. Then set
\ben
\bal
\la^{[1]}_m(u)&=(1+(\la_{m}^{(1)}+1)\tss u^{-1})\dots (1+(\la_{m}^{(k)}+1)\tss u^{-1})
(1+\la_{m}^{(k+1)}u^{-1})\dots (1+\la_{m}^{(p)}u^{-1}),\\[0.4em]
\la^{[1]}_{m+1}(u)&=(1+(\la_{m+1}^{(1)}+1)\tss u^{-1})\dots (1+(\la_{m+1}^{(k)}+1)\tss u^{-1})
(1+\la_{m}^{(k+1)}u^{-1})\dots (1+\la_{m}^{(p)}u^{-1}).
\eal
\een
Note that the series coincide with those given by \eqref{lammpo}.
Then apply the odd reflection to the pair $(\la^{[1]}_m(u),\la_{m+2}(u))$
and then continue as described in \eqref{oddlami} to conclude that the highest weight
of the module $L(\la(u))$ associated with the parity sequence $\se=1\dots 10\dots 010$
has the form
\ben
\big(\la_1(u),\dots,\la_{m-1}(u),\la^{[1]}_{m+1}(u),\dots,\la^{[1]}_{m+n-1}(u),
\la^{[n-1]}_{m}(u),\la_{m+n}(u),\la_{(m+n)'}(u)\big).
\een
Finally, we use a composition of embeddings \eqref{embedgen} as in
\cite[Cor.~3.2]{m:dt} to get the embedding
\beql{embospyang}
\X(\osp_{2|2})\hra \X(\osp_{2n|2m}).
\eeq
The cyclic span $\X(\osp_{2|2})\xi$ is a finite-dimensional highest weight
representation with the highest weight
$(\la^{[n-1]}_{m}(u),\la_{m+n}(u),\la_{(m+n)'}(u))$. Therefore
conditions \eqref{condosp}
follow from Corollary~\ref{cor:fdim}, as pointed out in Remark~\ref{rem:fdcom}\tss(ii).
\epf

\bcj\label{conj:fdim}
The conditions on the highest weight $\la(u)$ stated in Theorem~\ref{thm:necess}
are sufficient for the representation $L(\la(u))$ to be finite-dimensional.
\ecj

In the next section we will prove Conjecture~\ref{conj:fdim}
for representations with linear highest weights.

\section{Representations with linear highest weights}
\label{sec:rl}

Observe that by twisting any highest weight representation $L(\la(u))$ of
the extended Yangian $\X(\osp_{2n|2m})$
by a suitable automorphism \eqref{muf}, we can get a highest weight representation
whose components satisfy the condition
\beql{lalapr}
\la_{m+n}(u)\ts\la_{(m+n)'}(u+1)=1.
\eeq
We will now assume this condition so that the highest weight $\la(u)$ is determined
by the tuple \eqref{glnmt}. We will be concerned with those tuples
where all components are linear in $u^{-1}$. By a slight abuse of language we will
refer to such tuples as {\em linear highest weights}, even though the component
$\la_{(m+n)'}(u)$ does not have to be linear.

We will suppose that
\beql{linhwlaze}
\la_i(u)=1+\la_i\ts u^{-1},\qquad i=1,\dots,m+n,\quad\la_i\in\CC.
\eeq
Hence by \eqref{lalapr},
\beql{lamnpr}
\la_{(m+n)'}(u)=\frac{u-1}{u+\la_{m+n}-1}.
\eeq

For complex numbers $a$ and $b$ we will write $a\rar b$ or $b\lar a$ to mean
that $a-b\in\ZZ_+$.

\bth\label{thm:linear}
The representation $L(\la(u))$ of $\X(\osp_{2n|2m})$ with the highest weight
defined in \eqref{linhwlaze} and \eqref{lamnpr}
is finite-dimensional if and only if
\beql{ydiag}
\la_1\lar\cdots\lar\la_m=-l\fand \la_{m+1}\rar\cdots\rar\la_{m+l}\rar 0=\la_{m+l+1}=\cdots=\la_{m+n}
\eeq
for some $l\in\ZZ_+$; if $l\geqslant n$ then the second part of condition \eqref{ydiag}
is understood as
\ben
\la_{m+1}\rar\cdots\rar\la_{m+n}\rar 0.
\een
\eth

\bpf
Suppose first that $\dim L(\la(u))<\infty$. By Theorem~\ref{thm:necess}
the highest weight must satisfy conditions
\eqref{condgl}--\eqref{condosp}. Observe that the condition
\ben
1+au^{-1}\rar 1+bu^{-1},\qquad a,b\in\CC,
\een
for series in $u^{-1}$
is equivalent to $a\rar b$. Hence, we must have
\ben
\la_1\lar\cdots\lar\la_m\Fand \la_{m+1}\rar\cdots\rar\la_{m+n}.
\een
Furthermore, by condition \eqref{condo} we also have
\ben
1+\la_{m+n-1}\ts u^{-1}\rar \frac{u-1}{u+\la_{m+n}-1}
\een
and so
\ben
\frac{(u+\la_{m+n-1})(u+\la_{m+n}-1)}{u(u-1)}=\frac{Q(u+1)}{Q(u)}
\een
for a monic polynomial $Q(u)$ in $u$. This implies that $\la_{m+n}\rar 0$.

Now we make use of condition \eqref{condosp} and assume first that
the chain of odd reflections yields
\beql{aasge}
\la^{[n-1]}_{m}(u)=1+(\la_m+n-1)u^{-1}.
\eeq
Assuming further that $\la_m+n-1\ne \la_{m+n}$ we come to
the relation \eqref{fdco}
with $p=1$ for the parameters
\ben
\la_1(u):=\la^{[n-1]}_{m}(u),\qquad \la_2(u):=\la_{m+n}(u),\qquad \la_{2'}(u):=\la_{(m+n)'}(u),
\een
so that
\ben
\frac{u(u-1)}{(u+\la_{m}+n)(u+\la_{m}+n-1)}
=\frac{P(u+2)}{P(u)}
\een
for a monic polynomial $P(u)$ in $u$. This is possible only if
$\la_m+n\lar 0$ so that \eqref{ydiag} holds with $l\geqslant n$.

The above assumptions can only be violated if the equality $\la_m+l=\la_{m+l+1}$
occurs for some $l\in\{0,1,\dots,n-1\}$. In this case, $\la_m+l\in\ZZ_+$ and
\beql{aasgena}
\la^{[n-1]}_{m}(u)=1+(\la_m+l+s)u^{-1},
\eeq
for some $s\in\ZZ_+$ by the definition
of odd reflections in \eqref{oddrefl}. Hence condition \eqref{condosp} can only hold
if $s=0$ and
\ben
\la_m+l=\la_{m+l+1}=\dots=\la_{m+n}=0,
\een
implying \eqref{ydiag}. Note that the case $l=n-1$ is included by
\eqref{fdco}
with $p=0$ for the parameters
\ben
\la_1(u):=1,\qquad \la_2(u):=1,\qquad \la_{2'}(u):=\frac{u(u-1)}{(u+\la_{m+n})(u+\la_{m+n}-1)},
\een
implying $\la_{m+n}=0$.

Conversely, suppose that the conditions on the
components of the highest weight $\la(u)$ are satisfied.
We will interpret the highest weight in terms of
a Young diagram $\Ga=(\Ga_1,\Ga_2,\dots)$, contained
in the $(m,n)$-hook; cf. \cite[Sec.~2.1]{cw:dr}.
This means that the digram
$\Ga$ satisfies the
condition
$\Ga_{m+1}\leqslant n$, as illustrated:

\bigskip
\bigskip

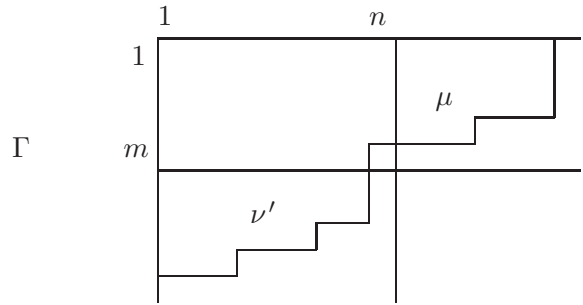
\begin{center}
\begin{picture}(150,90)
\thinlines

\put(0,0){\line(0,1){100}}
\put(90,0){\line(0,1){100}}
\put(0,50){\line(1,0){160}}
\put(0,100){\line(1,0){160}}

\put(0,10){\line(1,0){30}}
\put(30,10){\line(0,1){10}}
\put(30,20){\line(1,0){30}}
\put(60,20){\line(0,1){10}}
\put(60,30){\line(1,0){20}}
\put(80,30){\line(0,1){30}}
\put(80,60){\line(1,0){40}}
\put(120,60){\line(0,1){10}}
\put(120,70){\line(1,0){30}}
\put(150,70){\line(0,1){30}}

\put(0,105){\small$1$}
\put(80,105){\small$n$}
\put(-10,90){\small$1$}
\put(-13,55){\small$m$}

\put(-55,55){\small$\Ga$}

\put(105,75){\small$\mu$}
\put(35,30){\small$\nu\pr$}

\end{picture}
\end{center}


\noindent
The figure also shows partitions
$\mu=(\mu_1,\dots,\mu_m)$ and $\nu=(\nu_1,\dots,\nu_n)$
associated with $\Ga$. They are introduced by setting
\ben
\mu_i=\max\{\Ga_i-n,0\},\qquad i=1,\dots,m,
\een
and
\ben
\nu_j=\max\{\Ga\pr_j-m,0\},\qquad j=1,\dots,n,
\een
where $\Ga\pr$ denotes the conjugate partition so that $\Ga\pr_j$ is the length
of column $j$ in the diagram $\Ga$:

We will associate the $(m+n)$-tuple of integers $\Ga^{\sharp}$ with $\Ga$ by
\ben
\Ga^{\sharp}=(-\Ga_1,\dots,-\Ga_m,\nu_1,\dots,\nu_n).
\een
The conditions \eqref{ydiag} can now be stated equivalently as the relation
\beql{laga}
(\la_1,\dots,\la_m,\la_{m+1},\dots,\la_{m+n})=\Ga^{\sharp}
\eeq
for a suitable Young diagram $\Ga$, contained
in the $(m,n)$-hook. Hence, to complete the proof of the theorem,
it suffices to show that, given $\Ga$, the corresponding
representation $L(\la(u))$ of the extended Yangian $\X(\osp_{2n|2m})$ with
the parameters \eqref{laga} is finite-dimensional.
To this end, we use the same argument as in \cite[Sec.~4]{mr:rb}, so that
the representation $L(\la(u))$ is produced as a subquotient of the tensor product
of vector representations of the extended Yangian.

The {\em vector representation} of $\X(\osp^{\tss\se}_{2n|2m})$
(with an arbitrary parity sequence $\se$) on $\CC^{2n|2m}$ is defined by
\beql{vectre}
t_{ij}(u)\mapsto \de_{ij}+u^{-1}\tss
e_{ij}(-1)^{\bi}-(u+\ka)^{-1}\tss e_{j'i'}(-1)^{\bi\bj}\ts\ta_i\ta_j.
\eeq
The homomorphism property
follows from the $RTT$-relation \eqref{RTT} and
the Yang--Baxter equation satisfied by $R(u)$; cf. \cite{aacfr:rp}.

For the standard parity sequence $\se^{\st}$ and $d\in\{1,\dots,m\}$,
use the coproduct \eqref{Delta} to equip
the tensor product space $(\CC^{2n|2m})^{\ot d}$
with the action of $\X(\osp_{2n|2m})$ by setting
\ben
t_{ij}(u)\mapsto
\sum_{a_1,\dots,a_{d-1}=1}^{1'} t_{i,a_1}(u+d-1)\ot t_{a_1,a_2}(u+d-2)
\ot\dots\ot t_{a_{d-1},j}(u),
\een
where the generators act in the respective copies of the vector space
$(\CC^{2n|2m})^{\ot d}$ via the rule \eqref{vectre}.
Set
\beql{xid}
\xi_d=\sum_{\si\in\Sym_d} \sgn\si\cdot
e_{\si(1)}\ot\dots\ot e_{\si(d)} \in (\CC^{2n|2m})^{\ot d},
\eeq
where $\Sym_d$ denotes the symmetric group on $d$ symbols.
The calculations of \cite[Appendix]{m:rs} show that
the cyclic span $\X(\osp_{2n|2m})\ts \xi_d$ is a highest weight representation
of $\X(\osp_{2n|2m})$ with the highest weight $\la(u)$
whose components are found by
\ben
\la_i(u)=1-u^{-1}\quad\text{for}\quad i=1,\dots,d
\fand \la_i(u)=1\quad\text{for}\quad i=d+1,\dots,m+n+1.
\een
We will denote the irreducible quotient of this representation by $L^{\sharp\ts d}$.

We will use the following super-analogue of the well-known property
of Yangian modules, which is immediate from the coproduct formula \eqref{Delta}.

\ble\label{lem:cophw}
Let $L(\mu(u))$ and $L(\nu(u))$ be the irreducible highest weight $\X(\osp_{2n|2m})$-modules
with the highest weights
\ben
\mu(u)=\big(\mu_1(u),\dots,\mu_{m+n+1}(u)\big)\Fand
\nu(u)=\big(\nu_1(u),\dots,\nu_{m+n+1}(u)\big)
\een
and the respective highest vectors $\eta$ and $\ze$.
Then the cyclic span $\X(\osp_{2n|2m})(\eta\ot\ze)$ is a highest weight module
with the highest weight
\beql{mhw}
\big(\mu_1(u)\nu_1(u),\dots,\mu_{m+n+1}(u)\nu_{m+n+1}(u)\big).
\eeq
Hence, if the modules $L(\mu(u))$ and $L(\nu(u))$ are finite-dimensional,
then so is the irreducible highest weight module
with the highest weight \eqref{mhw}.
\qed
\ele

Given an $(m,n)$-hook partition $\Ga$, consider the tensor product module
\ben
L_{\Ga}=\bigotimes_{d=1}^m \Big(L^{\sharp\tss d}_{\Ga_d-1}\ot L^{\sharp\tss d}_{\Ga_d-2}\ot\dots\ot
L^{\sharp\tss d}_{\Ga_{d+1}}\Big),
\een
where $\Ga_{m+1}$ should be replaced by $0$. According to Lemma~\ref{lem:cophw}, the cyclic
$\X(\osp_{2n|2m})$-span of the tensor product of the highest weight vectors of
the modules $L^{\sharp\tss d}_{a}$ is a highest weight module with the highest weight
given by \eqref{linhwlaze} with
\ben
\la_i=-\Ga_i\quad\text{for}\quad i=1,\dots,m
\Fand \la_i=0\quad\text{for}\quad i=m+1,\dots,m+n.
\een

As the next step,
apply a sequence of odd reflections of type $A$ as defined in \eqref{oddrefl}
and an $\X(\osp_{2|2})$-type odd reflection introduced in Sec.~\ref{subsec:or}
to calculate
the highest weight of the irreducible quotient $\overline L_{\Ga}$ of this cyclic span
associated with the parity sequence $\se=0\dots01\dots1$. By Proposition~\ref{prop:odd22},
if the highest weight of the $\X(\osp^{10}_{2|2})$-module $L(\la(u))$
has the form
\ben
\la(u)=(1+(\la_m+n-1)u^{-1},1,1),\qquad \la_m+n\lar 0,
\een
then the $\X(\osp^{01}_{2|2})$-module $L(\la(u))^{\vs}$
is isomorphic to the highest weight module $L(\wt\la(u))$ with
\ben
\wt\la(u)=\Big(1+u^{-1},1+(\la_m+n)u^{-1},\frac{u-1}{u+\la_m+n-1}\Big).
\een
Therefore, the highest weight of the $\X(\osp^{\se}_{2n|2m})$-module $\overline L_{\Ga}$
is found by
\ben
\Big(1+\ol\Ga\pr_1u^{-1},\dots,1+\ol\Ga\pr_nu^{-1},
1+(\la_1+n)u^{-1},\dots,1+(\la_m+n)u^{-1},\frac{u-1}{u+\la_m+n-1}\big),
\een
in the case $\la_m+n\lar 0$ (so that $\ol\Ga\pr_i=m$ for all $i=1,\dots,n$), and
\ben
\big(1+\ol\Ga\pr_1u^{-1},\dots,1+\ol\Ga\pr_nu^{-1},
1+(\la_1+n)u^{-1},\dots,1+(\la_l+n)u^{-1},1\dots,1,1\big),
\een
in the case $\la_m=-l$ for some $l\in\{0,1,\dots,n-1\}$, where
$\ol\Ga=(\Ga_1,\dots,\Ga_m)$ is the Young diagram with $m$ rows.

Furthermore,
for $d\in\{1,\dots,n\}$
equip
the tensor product space $(\CC^{2n|2m})^{\ot d}$
with the action of $\X(\osp^{\tss\se}_{2n|2m})$ by setting
\ben
t_{ij}(u)\mapsto
\sum_{a_1,\dots,a_{d-1}=1}^{1'} t_{i,a_1}(u-d+1)\ot t_{a_1,a_2}(u-d+2)
\ot\dots\ot t_{a_{d-1},j}(u),
\een
where the generators act in the respective copies of the vector space
$(\CC^{2n|2m})^{\ot d}$ via the rule \eqref{vectre}.
The vector $\xi_d$ defined by the same formula \eqref{xid}
now generates a highest weight representation
of $\X(\osp^{\tss\se}_{2n|2m})$ with the highest weight $\la(u)$
whose components are found by
\ben
\la_i(u)=1+u^{-1}\quad\text{for}\quad i=1,\dots,d
\fand \la_i(u)=1\quad\text{for}\quad i=d+1,\dots,m+n+1.
\een
We will denote the irreducible quotient of this representation by $L^{\flat\tss d}$.
If the parameter $l$ in \eqref{ydiag} exceeds $n-1$, it should be understood as equal to $n$
in the argument below; in that case we set $\nu_{n+1}:=0$.
Consider
the tensor product module
\ben
\overline L_{\Ga}\ot\bigotimes_{d=1}^l \Big(L^{\flat\tss d}_{-\nu_d-m+1}
\ot L^{\flat\tss d}_{-\nu_d-m+2}\ot\dots\ot
L^{\flat\tss d}_{-\nu_{d+1}-m}\Big).
\een
By Lemma~\ref{lem:cophw}, the cyclic
$\X(\osp^{\tss\se}_{2n|2m})$-span of the tensor product of the highest weight vectors of
the tensor factors is a highest weight module with the highest weight
obtained from that of $\overline L_{\Ga}$ by replacing the first $n$ components by
$\ol\Ga\pr_i\mapsto \nu_i$ for $i=1,\dots,n$.
Hence the highest weight of the corresponding $\X(\osp_{2n|2m})$-module
is given by
\eqref{linhwlaze} and \eqref{lamnpr}, with the relation \eqref{laga}.
This completes the proof
of the sufficiency of the conditions \eqref{ydiag}.
\epf

Theorem~\ref{thm:linear} thus
confirms Conjecture~\ref{conj:fdim} in the case
of linear highest weights.

\appendix

\section{Isomorphism between the Yangians associated with
the Lie superalgebras $\osp_{2|2}$ and $\gl_{1|2}$}
\label{sec:ib}

We start by recalling the Gauss decomposition
of the generators matrix $T(u)$ associated with the Yangian $\Y(\gl_{n|m})$
or the extended Yangian $\X(\osp_{2n|2m})$. Denoting by $N$
the size of the matrix $T(u)$, write
\ben
T(u)=F(u)\ts H(u)\ts E(u),
\een
where $F(u)$, $H(u)$ and $E(u)$ are uniquely determined matrices of the form
\ben
F(u)=\begin{bmatrix}
1&0&\dots&0\ts\\
f_{21}(u)&1&\dots&0\\
\vdots&\vdots&\ddots&\vdots\\
f_{N1}(u)&f_{N2}(u)&\dots&1
\end{bmatrix},
\qquad
E(u)=\begin{bmatrix}
\ts1&e_{12}(u)&\dots&e_{1N}(u)\ts\\
\ts0&1&\dots&e_{2N}(u)\\
\vdots&\vdots&\ddots&\vdots\\
0&0&\dots&1
\end{bmatrix},
\een
and $H(u)=\diag\ts\big[h_1(u),\dots,h_N(u)\big]$. The entries
of the matrices $F(u)$, $H(u)$ and $E(u)$ are expressed
by universal quasideterminant formulas \cite{gr:dm}.

As in Sec.~\ref{sec:hw}, the Yangian $\Y(\gl^{\ts 011}_{1|2})$ associated
with the parity sequence $011$
is defined as the $\ZZ_2$-graded algebra
with generators $\bar t_{ij}^{\ts(r)}$
of parity $\bi+\bj\mod 2$, where $1\leqslant i,j\leqslant 3$ and $r=1,2,\dots$;
here $\bi$ takes the respective values $0,1,1$ for $i=1,2,3$.
The defining relations are given
in \eqref{gldefrel} for the generators combined into series
\eqref{btexpa}. Its presentation in terms of the Gaussian generators was given
in \cite{g:gd} so that the Yangian $\Y(\gl^{\ts 011}_{1|2})$ is generated
by the coefficients of the series $\bar h_1(u),\bar h_2(u),\bar h_3(u)$ together with
$\bar e_{12}(u),\bar e_{23}(u)$ and $\bar f_{21}(u),\bar f_{32}(u)$. We use
barred notation to distinguish the Gaussian generators of the Yangian $\Y(\gl^{\ts 011}_{1|2})$
from those of the extended orthosymplectic Yangian.

The extended Yangian $\X(\osp^{\ts 01}_{2|2})$ associated with the parity sequence $01$
is generated by the coefficients of the series $h_1(u),h_2(u),h_3(u)$ together with
$e_{12}(u),e_{23}(u)$ and $f_{21}(u),f_{32}(u)$. The presentation in terms of these
generators is given in \cite[Sec.~5.2.4]{ft:oy}.

\bpr\label{prop:isomgg}
The mapping
\begin{align}
\bar h_i(u)&\mapsto h_1(2u-1)\ts h_i(2u),\qquad i=1,2,3,
\label{hi}\\
\bar e_{12}(u)&\mapsto e_{12}(2u),\quad \bar e_{23}(u)\mapsto \frac12\ts e_{23}(2u),
\label{ei}\\[-1em]
\intertext{and}
\bar f_{21}(u)&\mapsto f_{21}(2u),\quad \bar f_{32}(u)\mapsto f_{32}(2u),
\label{fi}
\end{align}
defines an isomorphism $\phi:\Y(\gl^{\ts 011}_{1|2})\to \X(\osp^{\ts 01}_{2|2})$.
\epr

\bpf
Recall the tensor product decomposition
\ben
\Y(\gl^{\ts 011}_{1|2})=\ZY(\gl^{\ts 011}_{1|2})\ot \Y(\sll^{\ts 011}_{1|2})
\een
established in \cite[Prop.~3]{g:gd}, where $\ZY(\gl^{\ts 011}_{1|2})$
denotes the center of the Yangian $\Y(\gl^{\ts 011}_{1|2})$, which is generated by
the coefficients of the series
\ben
\be(u)=\bar h_1(u)^{-1}\ts\bar h_2(u)\ts\bar h_3(u+1).
\een
Its inverse coincides with the {\em quantum Berezinian} as shown in \cite[Sec.~7]{g:gd}.

On the other hand, we have the decomposition
\ben
\X(\osp^{\ts 01}_{2|2})=\ZX(\osp^{\ts 01}_{2|2})\ot\Y(\osp^{\ts 01}_{2|2})
\een
implied by \eqref{tensordecom}, where the center
$\ZX(\osp^{\ts 01}_{2|2})$ is generated by
the coefficients of the series
\ben
\si(u)=h_1(u)^{-1}\ts h_2(u)\ts h_1(u+1)\ts h_3(u+2)
\een
which coincides with $c_V(u+1)$ in the notation of \cite[Lem.~4.45]{ft:oy}.
Hence, we derive from \eqref{hi} that the mapping takes
$\be(u)$ to $\si(2u)$,
thus yielding an isomorphism $\ZY(\gl^{\ts 011}_{1|2})\to \ZX(\osp^{\ts 01}_{2|2})$.

Furthermore, the respective subalgebras
$\Y(\sll^{\ts 011}_{1|2})\subset \Y(\gl^{\ts 011}_{1|2})$ and
$\Y(\osp^{\ts 01}_{2|2})\subset \X(\osp^{\ts 01}_{2|2})$
are generated by the coefficients of the series
\ben
\bar\ka_1(u)=\bar h_1(u)^{-1}\ts\bar h_2(u),\qquad
\bar\ka_2(u)=\bar h_2(u+1/2)^{-1}\ts\bar h_3(u+1/2),
\een
\ben
\bar\xi_1^+(u)=\bar f_{21}(u),\quad
\bar\xi_1^-(u)=\bar e_{12}(u),\quad
\bar\xi_2^+(u)=\bar f_{32}(u+1/2),\quad
\bar\xi_2^-(u)=-\bar e_{23}(u+1/2)
\een
and
\ben
\ka_1(u)= h_1(u)^{-1}\ts h_2(u),\qquad
\ka_2(u)= h_2(u+1)^{-1}\ts h_3(u+1),
\een
\ben
\xi_1^+(u)= f_{21}(u),\quad
\xi_1^-(u)= e_{12}(u),\quad
\xi_2^+(u)= f_{32}(u+1),\quad
\xi_2^-(u)=-\frac12\ts e_{23}(u+1);
\een
see \cite[Sec.~9]{g:gd} and \cite[Sec.~6.2]{ft:oy}, respectively.
The presentations of the algebras $\Y(\sll^{\ts 011}_{1|2})$ and
$\Y(\osp^{\ts 01}_{2|2})$ given in {\em loc. cit.} show that
the mapping defined in the proposition provides an isomorphism
between these algebras given by
\beql{isomdr}
\Y(\sll^{\ts 011}_{1|2})\to \Y(\osp^{\ts 01}_{2|2}),\qquad
\bar\ka_i(u)\mapsto \ka_i(2u),\quad \bar\xi_i^{\pm}(u)\mapsto \xi_i^{\pm}(2u),\quad i=1,2,
\eeq
completing the proof.
\epf

We can find the images of the generators of the Yangian
$\Y(\gl^{\ts 011}_{1|2})$ under the isomorphism of Proposition~\ref{prop:isomgg}
in terms of the $RTT$ presentations.

\bco\label{cor:rttim}
Under the isomorphism $\phi$ of Proposition~\ref{prop:isomgg} we have
\ben
\bal
\bar t_{11}(u)&\mapsto t_{11}(2u-1)\ts t_{11}(2u),\qquad
\bar t_{12}(u)\mapsto t_{11}(2u-1)\ts t_{12}(2u),\\[0.4em]
\bar t_{21}(u)&\mapsto t_{21}(2u)\ts t_{11}(2u-1),\qquad
\bar t_{13}(u)\mapsto \frac12\ts t_{11}(2u-1)\ts t_{13}(2u),\\[0.4em]
\bar t_{31}(u)&\mapsto t_{31}(2u)\ts t_{11}(2u-1),
\eal
\een
together with
\ben
\bar t_{i2}(u)\mapsto t_{i2}(2u)\ts t_{11}(2u-1)+t_{i1}(2u)\ts t_{12}(2u-1)
-t_{i1}(2u)\ts t_{11}(2u)^{-1}\ts t_{11}(2u-1)\ts t_{12}(2u)
\een
for $i=2,3$, and
\ben
\bar t_{23}(u)\mapsto \frac12\ts \big(t_{23}(2u)\ts t_{11}(2u-1)+t_{21}(2u)\ts t_{13}(2u-1)
-t_{21}(2u)\ts t_{11}(2u)^{-1}\ts t_{11}(2u-1)\ts t_{13}(2u)\big),
\een
while the image of $\bar t_{33}(u)$ is found from the formula
\begin{multline}
\phi(\bar t_{33}(u))-\sum_{i,j=1}^2 \phi(\bar t_{3i}(u))\tss \phi(\bar A_{ij}(u))
\ts \phi(\bar t_{j\tss 3}(u))\\
{}=t_{11}(2u-1)\Big(t_{33}(2u)-\sum_{i,j=1}^2 t_{3i}(2u)\tss A_{ij}(2u)\ts t_{j\tss 3}(2u)\Big),
\non
\end{multline}
where the matrix $A(u)=[A_{ij}(u)]$ and its barred counterpart are defined in \eqref{htta}.
\eco

\bpf
The relations are immediate from Proposition~\ref{prop:isomgg}
and the quasideterminant formulas for the
Gaussian generators given in \eqref{gauh} and \eqref{htta} together with
\ben
e_{12}(u)=t_{11}(u)^{-1}\ts t_{12}(u),\qquad e_{23}(u)=h_2(u)^{-1}\ts
\begin{vmatrix} t_{1\tss 1}(u)&t_{1\tss 3}(u)\\
                         t_{2\tss 1}(u)&\boxed{t_{2\tss 3}(u)}
           \end{vmatrix}
\een
and
\ben
f_{21}(u)=t_{21}(u)\ts t_{11}(u)^{-1},\qquad f_{32}(u)=
\begin{vmatrix} t_{1\tss 1}(u)&t_{1\tss 2}(u)\\
                         t_{3\tss 1}(u)&\boxed{t_{3\tss 2}(u)}
           \end{vmatrix}\ts h_2(u)^{-1}.
\een
We also use the observation that $[t_{11}(u),t_{11}(v)]=0$ in
$\X(\osp^{\ts 01}_{2|2})$ and the following consequence of \eqref{defrel}:
\ben
t_{11}(2u-1)\ts t_{ij}(2u)+t_{i1}(2u-1)\ts t_{1j}(2u)
=t_{ij}(2u)\ts t_{11}(2u-1)+t_{i1}(2u)\ts t_{1j}(2u-1)
\een
which holds for all $i,j\in\{2,3\}$.
\epf

\subsection{Another proof of Corollary~\ref{cor:fdim}}
\label{subsec:apc}

We will use the assumptions of Proposition~\ref{prop:odd22}. By applying the isomorphism
\eqref{isom22} we find that
in the $\X(\osp^{01}_{2|2})$-module $L(\la(u))^{\vs}$ we have
\ben
t_{21}(u)\ts\xi=0\Fand t_{ij}(u)\ts\xi=0
\een
for all pairs $1\leqslant i<j\leqslant 1'$ except for $(i,j)=(1,2), (2\pr,1')$, whereas
\ben
t_{ii}(u)\ts\xi=\la_{\si(i)}(u)\ts\xi
\een
for $i=1,2,2\pr,1'$ with $\si=(12)(2^{\tss\prime}1')$. Recall that
the series $\la_{1'}(u)$ is determined by the highest weight $\la(u)$ by relation \eqref{zeoo}.

We will also need the eigenvalue of the series $h_2(u)$ on the highest vector $\xi$.
To calculate it we use an alternative expression for this series obtained by writing
\ben
\bal
h_2(u)\ts h_1(u+1)&=\big(t_{22}(u)-t_{21}(u)\ts t_{11}(u)^{-1}\ts t_{12}(u)\big)\ts t_{11}(u+1)\\
{}&=t_{22}(u)\ts t_{11}(u+1)-t_{21}(u)\ts t_{12}(u+1),
\eal
\een
where we used the relation $t_{12}(u)\ts t_{11}(u+1)=t_{11}(u)\ts t_{12}(u+1)$
implied by \eqref{defrel}. As another consequence of \eqref{defrel}, we then get
\beql{hhte}
h_2(u)\ts h_1(u+1)=t_{22}(u+1)\ts t_{11}(u)+t_{12}(u+1)\ts t_{21}(u).
\eeq
Hence
\beql{htxi}
h_2(u)\ts\xi=\la_1(u+1)\ts\la_2(u)\ts\la_2(u+1)^{-1}\ts\xi.
\eeq

Now apply the isomorphism $\phi:\Y(\gl^{\ts 011}_{1|2})\to \X(\osp^{\ts 01}_{2|2})$
provided by Proposition~\ref{prop:isomgg} to regard $L(\la(u))^{\vs}$
as a $\Y(\gl^{\ts 011}_{1|2})$-module. By Corollary~\ref{cor:rttim},
in this module we have
\ben
\bar t_{21}(u)\ts\xi=0,\qquad \bar t_{23}(u)\ts\xi=0,\qquad \bar t_{13}(u)\ts\xi=0,
\een
and
\ben
\bar t_{11}(u)\ts\xi=\la_2(2u-1)\ts\la_2(2u)\ts\xi,\qquad
\bar t_{33}(u)\tss\xi=\la_2(2u-1)\ts\la_{1'}(2u)\ts\xi.
\een
To determine the action of $\bar t_{22}(u)$ on $\xi$, we use relation \eqref{hhte}
which holds in the same form for the respective barred series $\bar h_i(u)$ and $\bar t_{ij}(u)$
in $\Y(\gl^{\ts 011}_{1|2})$. We have
\ben
\la_2(2u)\ts\la_2(2u-1)\ts\bar t_{22}(u+1)\ts\xi=\la_2(2u+2)\ts\la_2(2u+1)\ts\bar h_{2}(u)\ts\xi.
\een
However, the image of $\bar h_{2}(u)$ under the isomorphism $\phi$
is found as $h_1(2u-1)\ts h_2(2u)$ by \eqref{hi} so that using \eqref{htxi} we get
\ben
\bar h_{2}(u)\ts\xi=\la_1(2u+1)\ts\la_2(2u-1)\ts\la_2(2u)\ts\la_2(2u+1)^{-1}\ts\xi,
\een
thus yielding
\ben
\bar t_{22}(u)\ts\xi=\la_1(2u-1)\ts\la_2(2u)\ts\xi.
\een

As the next step, apply the isomorphism
\ben
\Y(\gl^{\ts 011}_{1|2})\to \Y(\gl^{\ts 101}_{1|2}),\qquad \bar t_{ij}(u)\mapsto \bar t_{\si(i),\si(j)}(u),
\een
where $\si=(1\ts 2)$ is the transposition. The same vector space $L(\la(u))$ can now
be regarded as the highest weight $\Y(\gl^{\ts 101}_{1|2})$-module
whose highest weight is the triple of formal series
\ben
\big(\la_1(2u-1)\ts\la_2(2u),\la_2(2u-1)\ts\la_2(2u),\la_2(2u-1)\ts\la_{1'}(2u)\big).
\een
Take its composition with the automorphism of $\Y(\gl^{\ts 101}_{1|2})$ defined by
\ben
\bar t_{ij}(u)\mapsto f(u)\ts \bar t_{ij}(u)\qquad\text{with}\quad
f(u)=\Big(\frac{2u-1}{2u}\Big)^p\ts\la_2(2u)^{-1}.
\een
The highest weight of the resulting representation
will be
\ben
\Bigg(\Big(\frac{2u-1}{2u}\Big)^p\ts\la_1(2u-1),\Big(\frac{2u-1}{2u}\Big)^p\ts\la_2(2u-1),
\Big(\frac{2u-1}{2u}\Big)^p\ts\frac{\la_{2}(2u-1)\tss\la_{1'}(2u)}{\la_2(2u)}\Bigg).
\een
Under the assumptions of Proposition~\ref{prop:odd22}, the application
of the odd reflection of \cite[Cor.~4.6]{m:or} yields
the highest weight $\Y(\gl^{\ts 011}_{1|2})$-module
whose highest weight is the triple
\ben
\Bigg(\Big(\frac{2u+1}{2u}\Big)^p\ts\la_2(2u+1),\Big(\frac{2u+1}{2u}\Big)^p\ts\la_1(2u+1),
\Big(\frac{2u-1}{2u}\Big)^p\ts\frac{\la_{2}(2u-1)\tss\la_{1'}(2u)}{\la_2(2u)}\Bigg).
\een
By the results of \cite{zh:sy} (also reproduced in \cite[Thm~4.4]{m:or}), this
representation if finite-dimensional if and only if the ratio of the last two components
\ben
\Big(\frac{2u-1}{2u+1}\Big)^p\ts\frac{\la_{2}(2u-1)\tss\la_{1'}(2u)}{\la_2(2u)\ts\la_1(2u+1)}
\een
coincides with $Q(u+1)/Q(u)$ for a monic polynomial $Q(u)$ in $u$. By replacing $u$ by $u/2$
and using \eqref{zeoo} we come to the condition \eqref{fdco}, thus completing the proof.

\section*{Declarations}

\subsection*{Competing interests}
The author has no competing interests to declare that are relevant to the content of this article.

\subsection*{Funding}
This work was supported by the Australian Research Council, grant DP240101572.

\subsection*{Availability of data and materials}
No data was used for the research described in the article.

\bigskip
\bigskip

\small
\noindent
School of Mathematics and Statistics\newline
University of Sydney,
NSW 2006, Australia\newline
alexander.molev@sydney.edu.au

\end{document}